\documentstyle{amsppt}
\voffset-13mm
\magnification1200
\pagewidth{130mm}
\pageheight{208mm}
\hfuzz=2.5pt\rightskip=0pt plus1pt
\binoppenalty=10000\relpenalty=10000\relax
\TagsOnRight
\nologo
\addto\eightpoint{\normalbaselineskip=.92\normalbaselineskip\normalbaselines}
\let\[\lfloor
\let\]\rfloor
\let\kappa\varkappa
\let\phi\varphi
\let\ol\overline
\redefine\d{\roman d}
\redefine\ord{\operatorname{ord}}
\define\Res{\operatorname{Res}}
\topmatter
\title
Irrationality of values of zeta-function
\endtitle
\author
W.~Zudilin
\endauthor
\endtopmatter
\rightheadtext{Irrationality of zeta-values}
\footnote""{2000 {\it Mathematics Subject Classification}.\enspace
Primary 11J72; Secondary 33C60.}

\subhead\indent
1. Introduction
\endsubhead
The irrationality of values of the zeta-function~$\zeta(s)$
at odd integers $s\ge3$ is one of the most attractive problems
in number theory. Inspite of a deceptive simplicity
and more than two-hundred-year history of the problem,
all done in this direction can easily be counted.
It was only~1978, when Ap\'ery~\cite{A} obtained the irrationality
of~$\zeta(3)$ by a presentation of ``nice'' rational
approximations to this number. During next years
the phenomenon of Ap\'ery's sequence was recomprehended
more than once from positions of different analytic methods
(see~\cite{N2} and the bibliography cited there);
new approaches gave rise to improve Ap\'ery's result
{\it quantitatively}, i.e., to get a ``sharp'' irrationality
measure of~$\zeta(3)$ (last stages in this direction
are the articles~\cite{H2},~\cite{RV}).
Finally, in~2000 Rivoal~\cite{R1} constructed
linear forms with rational coefficients involving
values of~$\zeta(s)$ only at odd integers $s>1$
and proved that {\it there exist infinitely many
irrational numbers among
$\zeta(3),\zeta(5),\zeta(7),\dots$\rom;
more precisely\rom, for the dimension $\delta(a)$
of spaces spanned over~$\Bbb Q$ by
$1,\zeta(3),\zeta(5),\dots,\zeta(a-2),\zeta(a)$\rom,
where $a$~is odd\rom, there holds the estimate\/}
$$
\delta(a)
\ge\frac{\log a}{1+\log2}\bigl(1+o(1)\bigr)
\qquad \text{\it as $a\to\infty$}.
$$

\subhead\indent
2. Main results
\endsubhead
In this note we generalize Rivoal's construction~\cite{R1}
and prove the following results.

\proclaim{\indent\smc Theorem 1}
Each of the following collections
$$
\gathered
\{\zeta(5), \; \zeta(7), \; \zeta(9), \; \zeta(11), \; \zeta(13), \;
\zeta(15), \; \zeta(17), \; \zeta(19), \; \zeta(21)\},
\\
\{\zeta(7), \; \zeta(9), \; 
\dots, \; \zeta(35), \; \zeta(37)\},
\qquad
\{\zeta(9), \; \zeta(11), \; 
\dots, \; \zeta(51), \; \zeta(53)\}
\endgathered
\tag1
$$
contains at least one irrational number.
\footnote"$^1$"{%
After finishing this paper the author knew that
Rivoal~\cite{R3} had independently obtained the claim of Theorem~1
for the first collection in~\thetag{1} by another
generalization of his construction from~\cite{R1}.
}
\endproclaim

\proclaim{\indent\smc Theorem 2}
For each odd integer $b\ge1$ the collection
$$
\zeta(b+2), \; \zeta(b+4), \; \dots, \; \zeta(8b-3), \; \zeta(8b-1)
$$
contains at least one irrational number.
\endproclaim

\proclaim{\indent\smc Theorem 3}
There exist odd integers $a_1\le145$ and $a_2\le1971$ such
that the numbers $1,\zeta(3),\zeta(a_1),\zeta(a_2)$
are linearly independent over~$\Bbb Q$.
\endproclaim

Theorem~3 improves corresponding result from~\cite{R2},
where the linear independence of numbers
$1,\zeta(3),\zeta(a)$ was established for some $a\le169$.

\proclaim{\indent\smc Theorem 4}
For each odd integer $a\ge3$ there holds the absolute estimate
$$
\delta(a)>0.395\,\log a>\frac23\cdot\frac{\log a}{1+\log2}.
\tag2
$$
\endproclaim

We stress that our proofs of Theorems~1--4 exploit calculations
via the saddle point method (Section~4)
and ideologically leans on the works~\cite{N2},~\cite{He}.
An improvelment of arithmetic estimates (i.e.,
of denominators of numerical linear forms)
in the spirit of~\cite{H2},~\cite{RV} (Section~3)
allows us to sharpen the lower estimate of $\delta(a)$
in Theorems~3,~4 for small values of~$a$.
In Section~5 we obtain not only an upper bound
but also precise asymptotics of coefficients
of linear forms. Finally, we prove Theorems~1--4 in Section~6.

The main results of the work were announced in the communication~\cite{Z}.

\medskip
The author is grateful to Professor Yu.\,V.~Nesterenko
for his permanent attention to the work.
This research was carried out with
the partial support of the INTAS--RFBR grant no.~IR-97-1904.

\subhead\indent
3. Analytic construction
\endsubhead
We fix positive odd parameters $a,b,c$ such that
$c\ge3$, $a>b(c-1)$,
and for each positive integer~$n$
consider the rational function
$$
\align
R(t)=R_n(t)
&:=\frac{\bigl((t\pm(n+1))\dotsb(t\pm cn)\bigr)^b}
{\bigl(t(t\pm1)\dotsb(t\pm n)\bigr)^a}\cdot(2n)!^{a+b-bc}
\\
&\phantom:=(-1)^n\cdot\biggl(\frac{\sin\pi t}\pi\biggr)^b
\cdot\frac{\Gamma(\pm t+cn+1)^b\Gamma(t-n)^{a+b}}{\Gamma(t+n+1)^{a+b}}
\cdot(2n)!^{a+b-bc},
\tag3
\endalign
$$
where the record~`$\pm$' means that the product contains factors
corresponding both to a sign~`$+$' and to a~`$-$'.
To the function~\thetag{3} assign the infinite sum
$$
I=I_n:=\sum_{t=n+1}^\infty\frac1{(b-1)!}\frac{\d^{b-1}R(t)}{\d t^{b-1}};
\tag4
$$
the series on the right-hand side of~\thetag{4} converges
absolutely since $R(t)=O(t^{-2})$ as $t\to\infty$.
Decomposing the function~\thetag{3} in a sum of partial fractions
and using its oddness, we deduce that
$$
I=\sum\Sb\text{$s$ is odd}\\b<s<a+b\endSb A_s\zeta(s)-A_0
\tag5
$$
(see~\thetag{10} below), where denominators
of the rational numbers $A_s=A_{s,n}$ grow not faster
than exponentially (cf.~\cite{R1}, lemmes~1,~5).
By~$D_n$ denote the least common multiple of numbers
$1,2,\dots,n$; the prime number theorem yields
$$
\lim_{n\to\infty}\frac{\log D_n}n=1.
$$

\proclaim{\indent\smc Lemma 1}
For each odd integer $c\ge3$ there exists a sequence
of integers $\Pi_n=\allowmathbreak\Pi_{n,c}\ge1$\rom,
$n=1,2,\dots$\rom, such
that the numbers $\Pi_n^{-b}D_{2n}^{a+b-1}A_{s,n}$
are integral and the limit relation
$$
\varpi_c:=\lim_{n\to\infty}\frac{\log\Pi_{n,c}}n
=-\sum_{l=1}^{(c-1)/2}\biggl(2\psi\biggl(\frac{2l}{c-1}\biggr)
+2\psi\biggl(\frac{2l}c\biggr)+\frac{2c-1}l\biggr)
+2(c-1)(1-\gamma)
\tag6
$$
holds\rom; here $\gamma\approx0{.}57712$ is Euler's constant
and $\psi(x)$~is the logarithmic derivative of the gamma-function.
\endproclaim

\demo{\indent\smc Proof}
Let
$$
\Pi_n=\prod_{\sqrt{(c+1)n}<p\le2n}p^{\nu_p},
\qquad\text{where}\quad
\nu_p=\min_{k=0,\pm1,\dots,\pm n}\biggl\{\ord_p
\frac{(cn+k)!\,(cn-k)!}{(n+k)!^c(n-k)!^c}\biggr\}.
\tag7
$$
Then for rational functions
$$
G(t)=G_n(t):=\frac{(t\pm(n+1))\dotsb(t\pm cn)}
{\bigl(t(t\pm1)\dotsb(t\pm n)\bigr)^{c-1}},
\qquad
H(t)=H_n(t):=\frac{(2n)!}{t(t\pm1)\dotsb(t\pm n)}
$$
we have the inclusions
$$
\gathered
\Pi_n^{-1}\cdot\frac{D_{2n}^j}{j!}\,\frac{\d^j}{\d t^j}
\bigl(G(t)(t+k)^{c-1}\bigr)\bigg|_{t=-k}\in\Bbb Z,
\quad
\frac{D_{2n}^j}{j!}\,\frac{\d^j}{\d t^j}
\bigl(H(t)(t+k)\bigr)\bigg|_{t=-k}\in\Bbb Z,
\\
k=0,\pm1,\dots,\pm n,
\quad
j=0,1,2,\dots
\endgathered
\tag8
$$
(a proof of the inclusions~\thetag{8} for the function~$G(t)$
needs a certain generalization of the arithmetic scheme
of Nikishin--Rivoal). Now, representing the initial
function~\thetag{3} in the form $R(t)=G(t)^bH(t)^{a+b-bc}$
and applying Leibniz's rule for the differentiation of a product,
by~\thetag{8} we obtain
$$
\gathered
\Pi_n^{-b}D_{2n}^jB_{k,j}\in\Bbb Z,
\qquad\text{where}\quad
B_{k,j}=\frac1{j!}\,\frac{\d^j}{\d t^j}
\bigl(R(t)(t+k)^a\bigr)\bigg|_{t=-k},
\\
k=0,\pm1,\dots,\pm n,
\quad j=0,1,\dots,a-1.
\endgathered
\tag9
$$
These relations yield the desired inclusions
$\Pi_n^{-b}D_{2n}^{a+b-1}A_s\in\Bbb Z$
since
$$
\aligned
A_s&=(-1)^{b-1}\binom{s-1}{b-1}\sum_{k=-n}^nB_{k,a+b-s-1},
\qquad \text{$s$ is odd}, \quad b<s<a+b,
\\
A_0&=(-1)^{b-1}\sum_{k=-n}^n\sum_{l=1}^{k+n}\biggl(
\binom{a+b-2}{b-1}\frac{B_{k,0}}{l^{a+b-1}}
+\dots
+\binom b{b-1}\frac{B_{k,a-2}}{l^{b+1}}
+\frac{B_{k,a-1}}{l^b}\biggr).
\endaligned
\tag10
$$

By~\thetag{7}, for each prime number $p>\sqrt{(c+1)n}$
there holds
$$
\nu_p=\min_{|k|\le n}\phi_c\biggl(\frac np,\frac kp\biggr),
$$
where the function
$\phi_c(x,y)=\[cx+y\]+\[cx-y\]-c\[x+y\]-c\[x-y\]$
is periodic (of period~$1$) with respect to each of its arguments,
and $\[\,\cdot\,\]$~is the integral part of a number.
Direct calculations show us that
$$
\gathered
\min_{y\in\Bbb R}\phi_c(x,y)=\cases
2l-2 &\text{if $\dsize x-\[x\]\in\biggl[\frac{l-1}{c-1},\frac lc\biggr)
\cup\biggl[\frac12+\frac{l-1}{c-1},\frac12+\frac lc\biggr)$},
\\
2l-1 &\text{if $\dsize x-\[x\]\in\biggl[\frac lc,\frac l{c-1}\biggr)
\cup\biggl[\frac12+\frac lc,\frac12+\frac l{c-1}\biggr)$},
\endcases
\\
l=1,2,\dots,\frac{c-1}2.
\endgathered
$$
Now, applying the number prime theorem and following
arguments from~\cite{Ch}, Theorem~4.3 and Section~6;
\cite{H1}, Lemma~3.2, we obtain the limit relation~\thetag{6}.
This completes the proof.
\enddemo

It can easily be checked that the value~$\varpi_c$ in~\thetag{6}
behaves itself like $2c(1-\gamma)+O(\log c)$ as $c\to\infty$.

\subhead\indent
4. Asymptotics of linear forms
\endsubhead
Consider the functions
$$
\cot_bz=\frac{(-1)^{b-1}}{(b-1)!}\,
\frac{\d^{b-1}\cot z}{\d z^{b-1}},
\qquad b=1,2,\dots\,.
$$
Decomposing $\pi\,\cot\pi t$ in a sum of partial
fractions we see that for each integer $b\ge1$
$$
\pi^b\,\cot_b\pi t=\frac1{(t-k)^b}+O(1)
\tag11
$$
in a neighbourhood of $t=k\in\Bbb Z$.

\proclaim{\indent\smc Lemma 2}
For the value~\thetag{4} there holds the integral presentation
$$
I=-\frac1{2\pi i}\int_{M-i\infty}^{M+i\infty}
\pi^b\,\cot_b\pi t\cdot R(t)\,\d t,
\tag12
$$
where $M\in\Bbb R$~is an arbitrary positive constant
from the interval $n<M<cn$.
\endproclaim

\demo{\indent\smc Proof}
Consider the integrand in~\thetag{12} on a
rectangle~$\Cal P$ with vertices $M\pm iN$, $N+\frac12\pm iN$,
where an integer~$N$ is sufficiently large, $N>cn$.
Expanding the function~\thetag{3} in Taylor series in a neighbourhood
of~$t=k\in\Bbb Z$ and using the expansion~\thetag{11}
by Cauchy's theorem we obtain
$$
\frac1{2\pi i}\int_{\Cal P}\pi^b\,\cot_b\pi t\cdot R(t)\,\d t
=\sum_{M<k\le N}\Res_{t=k}\bigl(\pi^b\,\cot_b\pi t\cdot R(t)\bigr)
=\sum_{M<k\le N}\frac{R^{(b-1)}(k)}{(b-1)!}.
\tag13
$$
On the sides $[N+\frac12-iN,N+\frac12+iN]$,
$[M-iN,N+\frac12-iN]$, and $[N+\frac12+iN,M+iN]$
of the rectangle~$\Cal P$ there holds the relation $R(t)=O(N^{-2})$,
while the function $\cot_b\pi t$, which is a polynomial in~$\cot\pi t$,
is bounded.
Hence, tending~$N$ to~$\infty$ in~\thetag{13}
we get the desired presentation~\thetag{12}.
\enddemo

Our next claim follows from Lemma~2 after
the change of variables $t=n\tau$ and an application
of Stirling's formula to gamma-factors of the function~\thetag{3}.

\proclaim{\indent\smc Lemma 3}
For the sum~\thetag{4} there holds the asymptotic relation
$$
I=\tilde I
\cdot\frac{(-1)^n(2\sqrt{\pi n}\,)^{a+b-bc}(2\pi)^b}{n^{a-1}}
\cdot\bigl(1+O(n^{-1})\bigr)
\qquad \text{as $n\to\infty$},
$$
where
$$
\gather
\tilde I=\tilde I_n
:=-\frac1{2\pi i}\int_{\Cal M}
\sin^b\pi n\tau\cdot\cot_b\pi n\tau
\cdot e^{nf(\tau)}\cdot g(\tau)\,\d\tau,
\tag14
\\
\aligned
f(\tau)
&=(a+b-bc)2\log2
+b(\tau+c)\log(\tau+c)
+b(-\tau+c)\log(-\tau+c)
\\ \vspace{-1.5pt} &\qquad
+(a+b)(\tau-1)\log(\tau-1)
-(a+b)(\tau+1)\log(\tau+1),
\endaligned
\\
g(\tau)
=\frac{(\tau+c)^{b/2}(-\tau+c)^{b/2}}
{(\tau+1)^{(a+b)/2}(\tau-1)^{(a+b)/2}},
\endgather
$$
and the contour $\Cal M$~is a vertical line $\Re(\tau)=\mu$\rom,
$1<\mu<c$\rom, oriented from bottom to top.
\endproclaim

We mean the functions $f(\tau)$ and~$g(\tau)$
in the complex $\tau$-plane cut along
the rays $(-\infty,1]$ and $[c,+\infty)$, where
we choose that branches of the logarithm functions,
which we assume to take real values for $\tau\in(1,c)$.

For each $b\ge1$ the function $\sin^bz\cdot\cot_bz$
is a polynomial in~$\cos z$ with rational coefficients:
$$
\sin^bz\cdot\cot_bz=V_b(\cos z),
\qquad V_b(-y)=(-1)^bV_b(y), \quad \deg V_b=\max\{1,b-2\};
$$
this fact immediately follows from the relations
$$
V_1(y)=y, \qquad
V_{b+1}(y)=yV_b(y)+\frac1b(1-y^2)V_b'(y), \quad b=1,2,\dots\,.
$$
Consequently, the integral~\thetag{14} can be represented in the form
$$
\tilde I=-\sum\Sb k=-b\\\text{$k$ is odd}\endSb^bc_kJ_{n,k},
\tag15
$$
where $c_k=c_{-k}$ are some (rational) constants
satisfying $c_1=1$ for $b=1$, and $c_b=0$, $c_{b-2}\ne0$ for $b>1$;
$$
J_{n,\lambda}
=\frac1{2\pi i}\int_{\Cal M}e^{n(f(\tau)-\lambda\pi i\tau)}
\cdot g(\tau)\,\d\tau
=\ol{J_{n,-\lambda}},
\qquad -b\le\lambda\le b,
\tag16
$$
and the overline means the complex conjugation.
To calculate asymptotics of the integrals~\thetag{16},
we apply the saddle point method exchanging the contour of integration
$\Cal M:\Re(\tau)=\mu$ by a contour~$\Cal M_\lambda$,
which passes through a (unique) saddle point~$\tau_\lambda$
in such a way that the integrand
achieves its maximal value at~$\tau_\lambda$. Saddle points
in the domain $\Re(\tau)>0$ can be determined from the equation
$$
f'(\tau)=\lambda\pi i, \qquad \lambda\in\Bbb R.
\tag17
$$

It follows easily that the polynomial
$$
(\tau+c)^b(\tau-1)^{a+b}-(\tau-c)^b(\tau+1)^{a+b}
\tag18
$$
has at least one real root on the interval $(c,+\infty)$;
by~$\mu_1$ denote such the root nearest to $\tau=c$.

\proclaim{\indent\smc Lemma 4}
Suppose the root $\mu_1\in(c,+\infty)$ of the polynomial~\thetag{18}
satisfies the condition
$$
\mu_1
\le c+\frac{c^2-1}4\cdot\min\biggl\{\frac b{2(a+b)},\frac1{3c}\biggr\}.
\tag19
$$
Then all solutions of equation~\thetag{17}
in the domain $\Re(\tau)>0$ are exhausted by the
following list\rom:
\roster
\item"(a)" the ``real'' solution $\mu_1\pm i0$ for $\lambda=\pm b$\rom,
where a sign~`$+$' \rom(a sign~`$-$'\rom) in the record~$\pm i0$ coinsides
with the sign of~$\lambda$ and corresponds to upper
\rom(respectively\rom, lower\rom) bank of the cut $[c,+\infty)$\rom;
\item"(b)" a real solution $\mu_0\in(1,c)$ for $\lambda=0$\rom;
\item"(c)" a complex solution $\tau_\lambda\in(1,c)$ for
$0<|\lambda|<b$\rom; in addition\rom, the sign of~$\Im(\tau_\lambda)$
coinsides with the sign of~$\lambda$\rom,
and $\tau_\lambda=\ol{\tau_{-\lambda}}$.
\endroster
The set of solutions of~\thetag{17} generates a smooth closed curve
$$
\Re(f'(\tau))
=\log\frac{|\tau+c|^b|\tau-1|^{a+b}}{|\tau-c|^b|\tau+1|^{a+b}}
=0
\tag20
$$
in the domain $\Re(\tau)>0$\rom; this curve is contained between
two circles centered at $\tau=c$
of radii $(\mu_1-c)/2$\rom, $(c-\mu_0)/2$\rom;
there holds $\Re(f'(\tau))>0$ inside the curve
and $\Re(f'(\tau))<0$ outside it.
\endproclaim

\demo\nofrills{\indent\smc Proof} \
of the claim leans on a geometric interpretation of the function
$\Im(f'(\tau))$ and on a description of all solutions of~\thetag{17}
in the whole cut $\tau$-plane
(and not only in the domain $\Re(\tau)>0$).
\enddemo

With the use of Lemma~4 we choose the countour~$\Cal M_\lambda$
to calculate asymptotics of~$J_{n,\lambda}$ as $n\to\infty$
in the following way. If $\lambda>0$ then the countour~$\Cal M_\lambda$
consists of the vertical ray $(\mu_0-i\infty,\mu_0]$,
of the segment $[\mu_0,\mu_0+e^{i\theta}\sqrt{\mu_0^2-1}\,]$
passing through the saddle point~$\tau_\lambda$,
and of the horizontal ray
$[\mu_0+e^{i\theta}\sqrt{\mu_0^2-1},e^{i\theta}\sqrt{\mu_0^2-1}+\infty]$;
in the case $\lambda<0$ the contour~$\Cal M_\lambda$
is symmetric to~$\Cal M_{-\lambda}$ with respect
to real axis; lastly, the contour~$\Cal M_0$ remains
the vertical line $(\mu_0-i\infty,\mu_0+i\infty)$.
This choice of the contour~$\Cal M_\lambda$ and an application
of Laplace's method (see, e.g., \cite{B}, \S\,5.7)
yield the following claim.

\proclaim{\indent\smc Lemma 5}
Let $\lambda\in\Bbb R$\rom, $|\lambda|\le b$\rom,
and let $\tau_\lambda$~be the \rom(unique\rom)
solution of equation \thetag{17} in the domain $\Re(\tau)>0$.
Then there holds the asymptotic formula
$$
|J_{n,\lambda}|
=\frac{e^{n\Re(f_0(\tau_\lambda))}|g(\tau_\lambda)|}
{(2\pi n|f''(\tau_\lambda)|)^{1/2}}
\cdot\bigl(1+O(n^{-1})\bigr)
\qquad \text{as $n\to\infty$},
$$
where
$$
\align
f_0(\tau)
:=f(\tau)-f'(\tau)\tau
&=(a+b-bc)2\log2+bc\log(\tau+c)+bc\log(-\tau+c)
\\ &\qquad
-(a+b)\log(\tau+1)-(a+b)\log(\tau-1).
\endalign
$$
\endproclaim

\proclaim{\indent\smc Lemma 6}
Suppose condition~\thetag{19} is satisfied.
Then for the linear forms~\thetag{5} there holds the limit relation
$$
\kappa:=\limsup_{n\to\infty}\frac{\log|I_n|}n
=\Re(f_0(\mu))=\log\frac{2^{2(a+b-bc)}|\mu+c|^{bc}|\mu-c|^{bc}}
{|\mu+1|^{a+b}|\mu-1|^{a+b}},
\tag21
$$
where $\mu$~is a real root \rom(i.e.\rom, $\mu_1$\rom)
of the polynomial~\thetag{18}
from the interval $(c,+\infty)$ for $b=1$\rom, or a root of this polynomial
in the domain~$\Im(\tau)>0$ with a maximal possible part $\Re(\mu)$
for $b>1$. In the case $b=1$ the limit superior in~\thetag{21}
can be replaced by the ordinary one.
\endproclaim

\demo{\indent\smc Proof}
All solutions of~\thetag{17} in the domain $\Re(\tau)>0$
for odd values $\lambda=k$ are simultaneously roots
of the polynomial~\thetag{18}. A routine test shows
that condition~\thetag{19} provides the increase
of the function $\Re(f_0(\tau))$ viewed as a function of~$\Re(\tau)$
(or, equivalently, as a function of~$\lambda$)
on the curve~\thetag{20} in the domain $\Re(\tau)>0$,
$\Im(\tau)\ge0$; hence only asymptotics of~$J_{n,\pm1}$ if $b=1$
and of~$J_{n,\pm(b-2)}$ if $b>1$ influence on asymptotics
of the integral~\thetag{15}. The application of Lemmas~5 and~3
yields the desired relation~\thetag{21}.
\enddemo

\subhead\indent
5. Estimates for coefficients of linear forms
\endsubhead
The values~$B_{k,j}$, $k=0,\pm1,\allowmathbreak\dots,\pm n$,
$j=0,1,\dots,a-1$, defined in~\thetag{9}
satisfy inequalities
$$
|B_{k,j}|
\le\bigl(2(a+bc-b)n\bigr)^j\cdot\max_{k=0,\pm1,\dots,\pm n}|B_{k,0}|
=\bigl(2(a+bc-b)n\bigr)^j\cdot\frac{(cn)!^{2b}(2n)!^{a+b-bc}}{n!^{2(a+b)}}.
$$
Using relations~\thetag{10} and Stirling's formula,
we then get

\proclaim{\indent\smc Lemma 7}
For the coefficients $A_s=A_{s,n}$ of the linear forms~\thetag{5}
there holds the estimate
$$
\gathered
\varlimsup_{n\to\infty}\frac{\log|A_{s,n}|}n
\le2bc\log c+2(a+b-bc)\log2,
\\
\text{$s=0$ or $s=b+1,\dots,a+b-1$ is odd}.
\endgathered
$$
\endproclaim

It is not hard to prove that the integrals
$$
\frac1{\pi i}\int_{iM-\infty}^{iM+\infty}
\biggl(\frac{\sin\pi t}\pi\biggr)^kR(t)\,\d t,
\qquad k=2,4,6,\dots,a-1,
\tag22
$$
passing through a horizontal line $\Im(t)=M$ with arbitrary $M>0$,
are linear combinations of coefficients $A_{b+2},\dots,A_{a+b-1}$
of the forms~\thetag{5}. Therefore, an application of
asymptotics of the gamma-function in the domain
$\Im(t)\ge M_0>0$ (see~\cite{B}, \S\,6.5)
and of the saddle point method to the integrals~\thetag{22}
makes more precise (insignificantly) the estimate of Lemma~7.

\proclaim{\indent\smc Lemma 8}
Suppose the real root $\mu_1\in(c,+\infty)$
of the polynomial~\thetag{18} satisfies condition~\thetag{19}\rom,
and let $\eta\in(0,+i\infty)$~be
an imaginary root of this polynomial
with a minimal possible absolute value.
Then for the coefficients $A_s=A_{s,n}$ of
the linear forms~\thetag{5} there holds the estimate
$$
\gathered
\varlimsup_{n\to\infty}\frac{\log|A_{s,n}|}n
\le\Re(f_0(\eta))
=\log\frac{2^{2(a+b-bc)}|\eta+c|^{bc}|\eta-c|^{bc}}
{|\eta+1|^{a+b}|\eta-1|^{a+b}},
\\
\text{$s=0$ or $s=b+1,\dots,a+b-1$ is odd\/};
\endgathered
$$
moreover\rom, in the case $s=a+b-1$
the limit superior can be replaced by the ordinary one
and the inequality becomes the equality.
\endproclaim

\subhead\indent
6. Proofs of main results
\endsubhead
By Lemmas~1,~6, if $-b\varpi_c+2(a+b-1)+\varkappa<0$
then there exists at least one irrational number
among values of~$\zeta(s)$, where $s$~is odd and $b<s<a+b$.
Taking $a=19$, $b=3$, $c=3$;
$a=33$, $b=5$, $c=3$, and $a=47$, $b=7$, $c=3$
respectively for the collections in~\thetag{1},
we deduce Theorem~1.
In Theorem~2, to each odd integer $b\ge1$ we assign $a=7b$, $c=3$;
to conclude the proof, it remains to note that
a real root of~\thetag{18}
from the interval $(3,+\infty)$ coincides with the root
$\mu_1\approx3.02472$ of the polynomial
$(\tau+3)(\tau-1)^8-(\tau-3)(\tau+1)^8$,
and
$\kappa+2(a+b-1)-b\varpi_c
<\Re(f_0(\mu_1))+16b-b\varpi_3
<-0.047\cdot b<0$.

In the case $b=1$ the criterion of linear independence from~\cite{N1}
in the same way as in~\cite{R1} allows us to obtain the lower estimate
for the value~$\delta(a)$, i.e.,
$$
\delta(a)\ge1-\frac{\kappa(a,c)+2a-\varpi_c}
{2c\log c+2(a-c+1)\log2+2a-\varpi_c},
\tag23
$$
where $\kappa=\kappa(a,c)$ is defined in~\thetag{21}.
Taking $a=145$, $c=21$ and $a=1971$, $c=131$, by~\thetag{23}
we obtain the estimates $\delta(145)\ge3$, $\delta(1971)\ge4$; in addition,
$\delta(3)=2$ due to~\cite{A}. This proves both Theorem~3
and Theorem~4 for $a<24999$. Further, for odd integers $a\ge20737=12^4+1$
we show stronger than~\thetag{2} estimate
$\delta(a)>\log_{12}a$ or, equivalently,
$$
\delta(12^m+1)>m, \qquad m=4,5,6,\dots,
\tag24
$$
choosing $c=2\cdot\[a/(3m^2)\]+1$ for each $a=12^m+1$
in~\thetag{23}. The estimate~\thetag{24} for $m=4,5,6,7$
is verified by direct calculations; finally, for $m\ge8$ we
use a trivial evaluation of the right-hand side in~\thetag{23}.
This completes the proof of Theorem~4.


\bgroup
\nobreak\vskip12pt minus6pt\indent\eightpoint\smc
\hbox to70mm{\vbox{\hsize=70mm%
\leftline{Moscow Lomonosov State University}
\leftline{Department of Mechanics and Mathematics}
\leftline{Vorobiovy Gory, Moscow 119899 RUSSIA}
\leftline{{\it E-mail address\/}: {\tt wadim\@ips.ras.ru}}
}}\par
\egroup
\newpage
\font\cyrmss=wncyr7 \relax
\font\cyrms=wncyr8 \relax
\font\cyris=wncyi8 \relax
\font\cyrm=wncyr10 \relax
\font\cyrb=wncyb10 \relax
\font\cyri=wncyi10 \relax
\font\cyrc=wncysc10 \relax
\font\cyrbx=wncyb10 scaled1200 \relax
\let\le\leqslant
\let\ge\geqslant
\let\limsup\varlimsup
\redefine\Re{\operatorname{Re}}
\redefine\Im{\operatorname{Im}}
\define\ctg{\operatorname{ctg}}
\catcode`\@=11
\firstpage@true
\catcode`\@=\active
\topmatter
\leftline{\cyrm UDK 511.3}
\vskip5mm
\title\cyrbx
O\lowercase{b irracional}\char"7E\lowercase{nosti
znaqeni}\char"1A{} \lowercase{dzeta-funkcii}
\endtitle
\author\cyrb
V.\,V.~Zudilin
\endauthor
\endtopmatter
\leftheadtext{\cyrms V.\,V.~Zudilin}
\rightheadtext{\cyrms
Irracional\char"5Enost\char"5E{}
znaqeni\char"12{} dzeta-funkcii}

\subhead\cyrb\indent
1. Vvedenie
\endsubhead\cyrm
Problema irracional\char"7Enosti znaqeni\char"1A{}
dzeta-funkcii~$\zeta(s)$ v neqetnyh toqkah $s\ge3$
yavlyaet\hbox{s}ya odno\char"1A{} iz samyh prityagatel\char"7Enyh
v teorii qisel. Nesmotrya na obmanqivuyu prostotu i bolee qem
dvuhvekovuyu isto\-riyu, poluqennye v e1tom napravlenii
rezul\char"7Etaty mozhno peresqitat\char"7E{} na pal\char"7Ecah.
Lix\char"7E{} v~1978\,g\. Aperi~\cite{A} udalos\char"7E{}
ustanovit\char"7E{} irracional\char"7Enost\char"7E{} $\zeta(3)$,
pred\char"7Fyaviv posledovatel\char"7Enost\char"7E{}
``horoxih'' racional\char"7Enyh
priblizheni\char"1A{} dlya e1togo qisla.
V dal\char"7Ene\char"1Axem fenomen posledovatel\char"7Enosti
Aperi byl ne\-odnokratno pereosmyslen s toqki zreniya razliqnyh
analitiqeskih metodov (sm.~\cite{N2} i citirovannuyu tam bibliografiyu);
novye podhody pozvolili usilit\char"7E{} rezul\char"7Etat Aperi
{\cyri koliqestvenno\/} -- poluqit\char"7E{} ``horoxuyu'' meru
irracional\char"7Enosti qisla~$\zeta(3)$
(poslednie e1tapy sorevnovaniya v e1tom napravlenii
-- raboty~\cite{H2},~\cite{RV}).
Nakonec, v~2000\,g\. Rivoal\char"7E~\cite{R1}
postroil line\char"1Anye formy s racional\char"7Enymi
koe1fficientami, soderzhawie znaqeniya~$\zeta(s)$ tol\char"7Eko
v neqetnyh toqkah $s>1$, i dokazal, qto {\cyri sredi qisel
$\zeta(3),\zeta(5),\zeta(7),\dots$
imeet\-\hbox{s}ya beskoneqno mnogo irracional\char"7Enyh\rom;
bolee toqno\rom, dlya razmernosti $\delta(a)$
pro\-stranstv\rom, natyanutyh nad~$\Bbb Q$ na qisla
$1,\zeta(3),\zeta(5),\allowmathbreak\dots,\zeta(a-2),\zeta(a)$\rom,
gde $a$~ne\-qetno\rom, spravedliva ocenka}
$$
\delta(a)
\ge\frac{\log a}{1+\log2}\bigl(1+o(1)\bigr)
\qquad \text{\cyri pri $a\to\infty$}.
$$

\subhead\cyrb\indent
2. Osnovnye rezul\char"7Etaty
\endsubhead\cyrm
V nastoyawe\char"1A{} zametke my obobwaem konstrukciyu Rivoalya~\cite{R1}
i dokazyvaem sleduyuwie rezul\char"7Etaty.

\proclaim{\cyrc\indent Teorema 1}\cyri
V kazhdom qislovom nabore
$$
\gathered
\{\zeta(5), \; \zeta(7), \; \zeta(9), \; \zeta(11), \; \zeta(13), \;
\zeta(15), \; \zeta(17), \; \zeta(19), \; \zeta(21)\},
\\
\{\zeta(7), \; \zeta(9), \; 
\dots, \; \zeta(35), \; \zeta(37)\},
\qquad
\{\zeta(9), \; \zeta(11), \; 
\dots, \; \zeta(51), \; \zeta(53)\}
\endgathered
\tag1
$$
imeet\hbox{s}ya po kra\char"1Ane\char"1A{}
mere odno irracional\char"7Enoe qislo.
\footnote"$^1$"{\cyrms
Posle zaverxeniya raboty nad stat\char"7Ee\char"1A{} avtoru stalo izvestno,
qto Rivoal\char"7E~\cite{R3} nezavisimo poluqil utverzhdenie teoremy~1
dlya pervogo iz naborov v~\thetag{1}, ispol\char"7Ezuya
inoe obobwenie konstrukcii iz~\cite{R1}.
}
\endproclaim

\proclaim{\cyrc\indent Teorema 2}\cyri
Dlya kazhdogo neqetnogo $b\ge1$ sredi qisel
$$
\zeta(b+2), \; \zeta(b+4), \; \dots, \; \zeta(8b-3), \; \zeta(8b-1)
$$
imeet\hbox{s}ya po kra\char"1Ane\char"1A{}
mere odno irracional\char"7Enoe.
\endproclaim

\proclaim{\cyrc\indent Teorema 3}\cyri
Suwestvuyut neqetnye $a_1\le145$ i $a_2\le1971$ takie\rom,
qto qisla $1,\zeta(3),\zeta(a_1),\zeta(a_2)$
line\char"1Ano nezavisimy nad~$\Bbb Q$.
\endproclaim\cyrm

Teorema~3 usilivaet sootvet\hbox{s}tvuyuwi\char"1A{}
rezul\char"7Etat raboty~\cite{R2},
gde usta\-novlena line\char"1Anaya nezavisimost\char"7E{} qisel
$1,\zeta(3),\zeta(a)$ dlya nekotorogo neqetnogo $a\le169$.

\proclaim{\cyrc\indent Teorema 4}\cyri
Dlya kazhdogo neqetnogo $a\ge3$ spravedliva absolyutnaya ocenka
$$
\delta(a)>0.395\,\log a>\frac23\cdot\frac{\log a}{1+\log2}.
\tag2
$$
\endproclaim\cyrm

Otmetim, qto dokazatel\char"7Estvo teorem~1--4 ispol\char"7Ezuet
vyqislenie asimptotiki s pomow\char"7Eyu metoda perevala
(p.~4) i ide\char"1Ano opiraet\hbox{s}ya na raboty~\cite{N2}, \cite{He}.
Usoverxenstvovanie arifmetiqeskih ocenok (znamenatele\char"1A{}
qislovyh line\char"1Anyh form) v duhe~\cite{H2},~\cite{RV},
privodimoe v~p.~3, pozvolilo utoqnit\char"7E{} ocenku snizu dlya $\delta(a)$
v teoremah~3,~4 pri malyh znaqeni\char"1A~$a$. V~p.~5 my po\-luqaem
ne tol\char"7Eko ocenku sverhu, no i toqnuyu asimptotiku koe1fficientov
line\char"1Anyh form. Nakonec, v~p.~6 my dokazyvaem teoremy~1--4.

Osnovnye rezul\char"7Etaty e1to\char"1A{} raboty
anonsirovany v soobwenii~\cite{Z}.

\medskip
Avtor iskrenne blagodaren professoru Yu.\,V.~Nesterenko
za postoyannoe vnimanie k rabote.
Nastoyawaya rabota vypolnena pri qastiqno\char"1A{} podderzhke
fonda {\rm INTAS} i Rossi\char"1Askogo fonda
fundamental\char"7Enyh issledovani\char"1A{}
(grant \rm no.~IR-97-1904).

\subhead\cyrb\indent
3. Analitiqeskaya konstrukciya
\endsubhead\cyrm
Zafiksiruem polozhitel\char"7Enye neqetnye parametry $a,b,c$,
$c\ge3$ i $a>b(c-1)$, i dlya kazhdogo celogo polozhitel\char"7Enogo~$n$
rassmotrim racional\char"7Enuyu funkciyu
$$
\align
R(t)=R_n(t)
&:=\frac{\bigl((t\pm(n+1))\dotsb(t\pm cn)\bigr)^b}
{\bigl(t(t\pm1)\dotsb(t\pm n)\bigr)^a}\cdot(2n)!^{a+b-bc}
\\
&\phantom:=(-1)^n\cdot\biggl(\frac{\sin\pi t}\pi\biggr)^b
\cdot\frac{\Gamma(\pm t+cn+1)^b\Gamma(t-n)^{a+b}}{\Gamma(t+n+1)^{a+b}}
\cdot(2n)!^{a+b-bc},
\tag3
\endalign
$$
gde znak `$\pm$' oznaqaet, qto v proizvedenii
uqastvuyut mnozhiteli, otveqayuwie kak znaku~`$+$', tak i~`$-$'.
Postavim v sootvet\hbox{s}tvie funkcii~\thetag{3} beskoneqnuyu summu
$$
I=I_n:=\sum_{t=n+1}^\infty\frac1{(b-1)!}\frac{\d^{b-1}R(t)}{\d t^{b-1}};
\tag4
$$
ryad v pravo\char"1A{} qasti~\thetag{4} \hbox{s}hodit\hbox{s}ya absolyutno, poskol\char"7Eku
$R(t)=O(t^{-2})$ pri $t\to\infty$.
Predstavlyaya funkciyu~\thetag{3} v vide summy proste\char"1Axih drobe\char"1A{}
i ispol\char"7Ezuya ee neqetnost\char"7E, zaklyuqaem, qto
$$
I=\sum\Sb\text{\cyrmss $s$ neqetno}\\b<s<a+b\endSb A_s\zeta(s)-A_0
\tag5
$$
(sm.~\thetag{10} dalee), gde znamenateli
racional\char"7Enyh qisel $A_s=A_{s,n}$ rastut
ne bystree qem e1ksponencial\char"7Eno (sr.~\cite{R1, \cyrm lemmy~1,~5}).
Oboznaqim qerez~$D_n$ naimen\char"7Exee obwee kratnoe qisel $1,2,\dots,n$;
iz asimptotiqeskogo zakona raspredeleniya prostyh qisel
$$
\lim_{n\to\infty}\frac{\log D_n}n=1.
$$

\proclaim{\cyrc\indent Lemma 1}\cyri
Dlya kazhdogo neqetnogo $c\ge3$ suwestvuet
posledovatel\char"7Enost\char"7E{}
celyh $\Pi_n=\Pi_{n,c}\ge1$\rom, $n=1,2,\dots$\rom, takaya\rom,
qto qisla $\Pi_n^{-b}D_{2n}^{a+b-1}A_{s,n}$
yavlyayut\hbox{s}ya celymi i spravedlivo
predel\char"7Enoe sootnoxenie
$$
\varpi_c:=\lim_{n\to\infty}\frac{\log\Pi_{n,c}}n
=-\sum_{l=1}^{(c-1)/2}\biggl(2\psi\biggl(\frac{2l}{c-1}\biggr)
+2\psi\biggl(\frac{2l}c\biggr)+\frac{2c-1}l\biggr)
+2(c-1)(1-\gamma),
\tag6
$$
gde $\gamma\approx0{.}57712$ -- postoyannaya E1\char"1Alera\rom,
a $\psi(x)$~-- logarifmiqeskaya proizvodnaya gamma-funkcii.
\endproclaim

\demo{\cyrc\indent Dokazatel\char"7Estvo}\cyrm
Polozhim
$$
\Pi_n=\prod_{\sqrt{(c+1)n}<p\le2n}p^{\nu_p},
\qquad\text{\cyrm gde}\quad
\nu_p=\min_{k=0,\pm1,\dots,\pm n}\biggl\{\ord_p
\frac{(cn+k)!\,(cn-k)!}{(n+k)!^c(n-k)!^c}\biggr\}.
\tag7
$$
Togda dlya racional\char"7Enyh funkci\char"1A{}
$$
G(t)=G_n(t):=\frac{(t\pm(n+1))\dotsb(t\pm cn)}
{\bigl(t(t\pm1)\dotsb(t\pm n)\bigr)^{c-1}},
\qquad
H(t)=H_n(t):=\frac{(2n)!}{t(t\pm1)\dotsb(t\pm n)}
$$
spravedlivy vklyuqeniya
$$
\gathered
\Pi_n^{-1}\cdot\frac{D_{2n}^j}{j!}\,\frac{\d^j}{\d t^j}
\bigl(G(t)(t+k)^{c-1}\bigr)\bigg|_{t=-k}\in\Bbb Z,
\quad
\frac{D_{2n}^j}{j!}\,\frac{\d^j}{\d t^j}
\bigl(H(t)(t+k)\bigr)\bigg|_{t=-k}\in\Bbb Z,
\\
k=0,\pm1,\dots,\pm n,
\quad
j=0,1,2,\dots
\endgathered
\tag8
$$
(dokazatel\char"7Estvo vklyuqeni\char"1A~\thetag{8}
dlya funkcii~$G(t)$ ispol\char"7Ezuet obobwenie
ari\-fmetiqesko\char"1A{} \hbox{s}hemy Nikixina--Rivoalya).
Zapisyvaya i\hbox{s}hodnuyu funkciyu~\thetag{3} v vide
$R(t)=G(t)^bH(t)^{a+b-bc}$ i primenyaya pravilo
Le\char"1Abnica dlya differencirovaniya proizvedeniya,
soglasno vklyuqeniyam~\thetag{8} poluqaem
$$
\gathered
\Pi_n^{-b}D_{2n}^jB_{k,j}\in\Bbb Z,
\qquad\text{\cyrm gde}\quad
B_{k,j}=\frac1{j!}\,\frac{\d^j}{\d t^j}
\bigl(R(t)(t+k)^a\bigr)\bigg|_{t=-k},
\\
k=0,\pm1,\dots,\pm n,
\quad j=0,1,\dots,a-1.
\endgathered
\tag9
$$
E1to daet trebuemye vklyuqeniya $\Pi_n^{-b}D_{2n}^{a+b-1}A_s\in\Bbb Z$,
tak kak
$$
\aligned
A_s&=(-1)^{b-1}\binom{s-1}{b-1}\sum_{k=-n}^nB_{k,a+b-s-1},
\qquad \text{\cyrm $s$ neqetno}, \quad b<s<a+b,
\\
A_0&=(-1)^{b-1}\sum_{k=-n}^n\sum_{l=1}^{k+n}\biggl(
\binom{a+b-2}{b-1}\frac{B_{k,0}}{l^{a+b-1}}
+\dots
+\binom b{b-1}\frac{B_{k,a-2}}{l^{b+1}}
+\frac{B_{k,a-1}}{l^b}\biggr).
\endaligned
\tag10
$$

Dlya kazhdogo prostogo $p>\sqrt{(c+1)n}$
soglasno~\thetag{7} vypolneno
$$
\nu_p=\min_{|k|\le n}\phi_c\biggl(\frac np,\frac kp\biggr),
$$
gde funkciya
$\phi_c(x,y)=\[cx+y\]+\[cx-y\]-c\[x+y\]-c\[x-y\]$
periodiqna (s periodom~$1$) po kazhdomu argumentu,
$\[\,\cdot\,\]$~-- celaya qast\char"7E{} qisla.
Proverka pokazyvaet, qto
$$
\gathered
\min_{y\in\Bbb R}\phi_c(x,y)=\cases
2l-2, &\text{\cyrm esli
$\dsize x-\[x\]\in\biggl[\frac{l-1}{c-1},\frac lc\biggr)
\cup\biggl[\frac12+\frac{l-1}{c-1},\frac12+\frac lc\biggr)$},
\\
2l-1, &\text{\cyrm esli
$\dsize x-\[x\]\in\biggl[\frac lc,\frac l{c-1}\biggr)
\cup\biggl[\frac12+\frac lc,\frac12+\frac l{c-1}\biggr)$},
\endcases
\\
l=1,2,\dots,\frac{c-1}2;
\endgathered
$$
ot\hbox{s}yuda s pomow\char"7Eyu asimptotiqeskogo zakona raspredeleniya prostyh qisel
i rassuzhdeni\char"1A{} iz~\cite{Ch, \cyrm teorema~4.3 i~\S\,6}, \cite{H1, \cyrm lemma~3.2}
poluqaem predel\char"7Enoe sootnoxenie~\thetag{6}. Lemma dokazana.
\enddemo

Neslozhno ubedit\char"7Esya v tom, qto veliqina~$\varpi_c$ v~\thetag{6}
pri $c\to\infty$ imeet poryadok $2c(1-\gamma)+O(\log c)$.

\subhead\cyrb\indent
4. Asimptotika line\char"1Anyh form
\endsubhead\cyrm
Opredelim funkcii
$$
\ctg_bz=\frac{(-1)^{b-1}}{(b-1)!}\,
\frac{\d^{b-1}\ctg z}{\d z^{b-1}},
\qquad b=1,2,\dots\,.
$$
Razlozhenie $\pi\,\ctg\pi t$ v summu proste\char"1Axih drobe\char"1A{}
pokazyvaet, qto dlya lyubogo celogo $b\ge1$ v okrestnosti
toqki $t=k\in\Bbb Z$ spravedlivo predstavlenie
$$
\pi^b\,\ctg_b\pi t=\frac1{(t-k)^b}+O(1).
\tag11
$$

\proclaim{\cyrc\indent Lemma 2}\cyri
Dlya veliqiny~\thetag{4} spravedlivo
integral\char"7Enoe predstavlenie
$$
I=-\frac1{2\pi i}\int_{M-i\infty}^{M+i\infty}
\pi^b\,\ctg_b\pi t\cdot R(t)\,\d t,
\tag12
$$
gde $M\in\Bbb R$~-- proizvol\char"7Enaya postoyannaya
iz intervala $n<M<cn$.
\endproclaim

\demo{\cyrc\indent Dokazatel\char"7Estvo}\cyrm
Rassmotrim podyntegral\char"7Enuyu funkciyu v~\thetag{12} na konture
pryamougol\char"7Enika~$\Cal P$ s verxinami $M\pm iN$, $N+\frac12\pm iN$,
gde celoe qislo~$N$ dostatoqno veliko, $N>cn$.
Raskladyvaya funkciyu~\thetag{3} v ryad Te\char"1Alora v ok\-restnosti
$t=k\in\Bbb Z$ i pol\char"7Ezuyas\char"7E{} predstavleniem~\thetag{11},
soglasno integ\-ral\char"7Eno\char"1A{} teoreme Koxi poluqaem
$$
\frac1{2\pi i}\int_{\Cal P}\pi^b\,\ctg_b\pi t\cdot R(t)\,\d t
=\sum_{M<k\le N}\Res_{t=k}\bigl(\pi^b\,\ctg_b\pi t\cdot R(t)\bigr)
=\sum_{M<k\le N}\frac{R^{(b-1)}(k)}{(b-1)!}.
\tag13
$$
Na storonah $[N+\frac12-iN,N+\frac12+iN]$,
$[M-iN,N+\frac12-iN]$ i $[N+\frac12+iN,M+iN]$
pryamougol\char"7Enika~$\Cal P$ vypolneno $R(t)=O(N^{-2})$ i funkciya
$\ctg_b\pi t$, yavlyayuwayasya mnogoqlenom ot~$\ctg\pi t$, ograniqena.
Poe1tomu predel\char"7Eny\char"1A{} perehod $N\to\infty$ v~\thetag{13}
privodit k~\thetag{12}.
\enddemo

Sleduyuwee utverzhdenie poluqaet\hbox{s}ya iz lemmy~2
posle zameny $t=n\tau$ i primeneniya k gamma-mnozhitelyam
funkcii~\thetag{3} formuly Stirlinga.

\proclaim{\cyrc\indent Lemma 3}\cyri
Pri $n\to\infty$ dlya summy~\thetag{4} vypolneno
$$
I=\tilde I
\cdot\frac{(-1)^n(2\sqrt{\pi n}\,)^{a+b-bc}(2\pi)^b}{n^{a-1}}
\cdot\bigl(1+O(n^{-1})\bigr),
$$
gde
$$
\gather
\tilde I=\tilde I_n
:=-\frac1{2\pi i}\int_{\Cal M}
\sin^b\pi n\tau\cdot\ctg_b\pi n\tau
\cdot e^{nf(\tau)}\cdot g(\tau)\,\d\tau,
\tag14
\\
\aligned
f(\tau)
&=(a+b-bc)2\log2
+b(\tau+c)\log(\tau+c)
+b(-\tau+c)\log(-\tau+c)
\\ \vspace{-1.5pt} &\qquad
+(a+b)(\tau-1)\log(\tau-1)
-(a+b)(\tau+1)\log(\tau+1),
\endaligned
\\
g(\tau)
=\frac{(\tau+c)^{b/2}(-\tau+c)^{b/2}}
{(\tau+1)^{(a+b)/2}(\tau-1)^{(a+b)/2}},
\endgather
$$
a kontur $\Cal M$~-- vertikal\char"7Enaya pryamaya $\Re\tau=\mu$\rom,
$1<\mu<c$\rom, prohodimaya snizu vverh.
\endproclaim\cyrm

My rassmatrivaem funkcii $f(\tau)$ i~$g(\tau)$
v $\tau$-ploskosti s razrezami vdol\char"7E{} luqe\char"1A{} $(-\infty,1]$
i $[c,+\infty)$, fiksiruya vetvi logarifmov,
prinimayuwie de\char"1Astvi\-tel\char"7Enye znaqeniya na intervale
$(1,c)$ vewestvenno\char"1A{} osi.

Dlya kazhdogo $b\ge1$ funkciya $\sin^bz\cdot\ctg_bz$
yavlyaet\hbox{s}ya mnogoqlenom ot~$\cos z$
s racional\char"7Enymi koe1fficientami:
$$
\sin^bz\cdot\ctg_bz=V_b(\cos z),
\qquad V_b(-y)=(-1)^bV_b(y), \quad \deg V_b=\max\{1,b-2\};
$$
e1tot fakt sleduet iz sootnoxeni\char"1A{}
$$
V_1(y)=y, \qquad
V_{b+1}(y)=yV_b(y)+\frac1b(1-y^2)V_b'(y), \quad b=1,2,\dots\,.
$$
Poe1tomu integral~\thetag{14} mozhno predstavit\char"7E{} v vide
$$
\tilde I=-\sum\Sb k=-b\\\text{\cyrmss $k$ neqetno}\endSb^bc_kJ_{n,k},
\tag15
$$
gde $c_k=c_{-k}$ -- nekotorye (racional\char"7Enye) postoyannye,
priqem $c_1=1$ dlya $b=1$ i $c_b=0$, $c_{b-2}\ne0$ dlya $b>1$,
$$
J_{n,\lambda}
=\frac1{2\pi i}\int_{\Cal M}e^{n(f(\tau)-\lambda\pi i\tau)}
\cdot g(\tau)\,\d\tau
=\ol{J_{n,-\lambda}},
\qquad -b\le\lambda\le b,
\tag16
$$
qerta sverhu oznaqaet kompleksnoe sopryazhenie. Dlya vyqisleniya
asimptotiki integralov~\thetag{16} my vospol\char"7Ezuemsya metodom
perevala, zamenyaya kontur integrirovaniya $\Cal M:\Re\tau=\mu$
na kontur~$\Cal M_\lambda$, prohodyawi\char"1A{} qerez (edinstvennuyu) toqku
perevala~$\tau_\lambda$, v kotoro\char"1A{} podyntegral\char"7Enaya funkciya
prinimaet maksimal\char"7Enoe znaqenie. Toqki perevala
v oblasti $\Re\tau>0$ udovletvoryayut uravneniyu
$$
f'(\tau)=\lambda\pi i, \qquad \lambda\in\Bbb R.
\tag17
$$

Kak neslozhno zametit\char"7E, mnogoqlen
$$
(\tau+c)^b(\tau-1)^{a+b}-(\tau-c)^b(\tau+1)^{a+b}
\tag18
$$
imeet po kra\char"1Ane\char"1A{} mere odin vewestvenny\char"1A{} koren\char"7E{}
na intervale $(c,+\infty)$; oboznaqim qerez~$\mu_1$
blizha\char"1Axi\char"1A{} iz e1tih korne\char"1A{} k toqke $\tau=c$.

\proclaim{\cyrc\indent Lemma 4}\cyri
Pust\char"7E{} koren\char"7E{}
$\mu_1\in(c,+\infty)$ mnogoqlena~\thetag{18}
udovletvoryaet ne\-ravenstvu
$$
\mu_1
\le c+\frac{c^2-1}4\cdot\min\biggl\{\frac b{2(a+b)},\frac1{3c}\biggr\}.
\tag19
$$
Togda v oblasti $\Re\tau>0$ vse rexeniya uravneniya~\thetag{17}
isqerpyvayut\hbox{s}ya sleduyuwim spiskom\rom:
\roster
\item"\cyrm a)"
``vewestvennoe'' rexenie $\mu_1\pm i0$ dlya $\lambda=\pm b$\rom,
gde znak~`$+$' \rom(znak~`$-$'\rom) v zapisi~$\pm i0$ sovpadaet so
zn\it\'a\cyri kom~$\lambda$ i otveqaet verhnemu
\rom(nizhnemu\rom) beregu razreza $[c,+\infty)$\rom;
\item"\cyrm b)"
vewestvennoe rexenie $\mu_0\in(1,c)$ dlya $\lambda=0$\rom;
\item"\cyrm v)"
kompleksnoe rexenie $\tau_\lambda\in(1,c)$ dlya
$0<|\lambda|<b$\rom, pri e1tom znak~$\Im\tau_\lambda$
sovpadaet so zn\it\'a\cyri kom~$\lambda$
i $\tau_\lambda=\ol{\tau_{-\lambda}}$.
\endroster
Mnozhestvo rexeni\char"1A{} uravneniya~\thetag{17}
obrazuyut v poluploskosti
$\Re\tau>0$ gladkuyu zamknutuyu krivuyu
$$
\Re f'(\tau)
=\log\frac{|\tau+c|^b|\tau-1|^{a+b}}{|\tau-c|^b|\tau+1|^{a+b}}
=0,
\tag20
$$
zaklyuqennuyu vnutri okruzhnoste\char"1A{} s centrom v toqke $\tau=c$
i radiusami $(\mu_1-c)/2$\rom, $(c-\mu_0)/2$\rom;
vnutri e1to\char"1A{} krivo\char"1A{} $\Re f'(\tau)>0$
i vne ee $\Re f'(\tau)<0$.
\endproclaim

\demo\nofrills{\cyrc\indent Dokazatel\char"7Estvo}\cyrm\
e1togo utverzhdeniya opiraet\hbox{s}ya na geometriqeskuyu interpretaciyu
funkcii $\Im f'(\tau)$ i opisanie rexeni\char"1A{} uravneniya~\thetag{17}
vo vse\char"1A{} $\tau$-plos\-kosti (a ne tol\char"7Eko v oblasti $\Re\tau>0$).
\enddemo

S pomow\char"7Eyu lemmy~4 dlya vyqisleniya asimptotiki~$J_{n,\lambda}$
pri $n\to\infty$ my vybiraem kontur~$\Cal M_\lambda$, sostoyawi\char"1A{}
v sluqae $\lambda>0$ iz vertikal\char"7Enogo luqa $(\mu_0-i\infty,\mu_0]$,
otrezka $[\mu_0,\mu_0+e^{i\theta}\sqrt{\mu_0^2-1}\,]$, prohodyawego
qerez toqku perevala $\tau_\lambda$, i gorizontal\char"7Enogo luqa
$[\mu_0+e^{i\theta}\sqrt{\mu_0^2-1},e^{i\theta}\sqrt{\mu_0^2-1}+\infty]$;
v sluqae $\lambda<0$ kontur~$\Cal M_\lambda$ simmetriqen~$\Cal M_{-\lambda}$
otnositel\char"7Eno vewestvenno\char"1A{} osi; nakonec, kontur~$\Cal M_0$
est\char"7E{} vertikal\char"7Enaya pryamaya $(\mu_0-i\infty,\mu_0+i\infty)$.
Tako\char"1A{} vybor kontura~$\Cal M_\lambda$ i primenenie metoda Laplasa
(sm., naprimer, \cite{B, \S\,5.7}) privodit k sleduyuwemu utverzhdeniyu.

\proclaim{\cyrc\indent Lemma 5}\cyri
Pust\char"7E{} $\lambda\in\Bbb R$\rom, $|\lambda|\le b$\rom,
i $\tau_\lambda$~-- \rom(edinstvennoe\rom) rexenie
uravneniya~\thetag{17} v oblasti $\Re\tau>0$.
Togda pri $n\to\infty$ spravedliva
asimptotiqeskaya formula
$$
|J_{n,\lambda}|
=\frac{e^{n\Re f_0(\tau_\lambda)}|g(\tau_\lambda)|}
{(2\pi n|f''(\tau_\lambda)|)^{1/2}}
\cdot\bigl(1+O(n^{-1})\bigr),
$$
gde
$$
\align
f_0(\tau)
:=f(\tau)-f'(\tau)\tau
&=(a+b-bc)2\log2+bc\log(\tau+c)+bc\log(-\tau+c)
\\ &\qquad
-(a+b)\log(\tau+1)-(a+b)\log(\tau-1).
\endalign
$$
\endproclaim

\proclaim{\cyrc\indent Lemma 6}\cyri
Pust\char"7E{} vypolneno uslovie~\thetag{19}.
Togda dlya line\char"1Anyh form~\thetag{5}
spravedlivo predel\char"7Enoe sootnoxenie
$$
\kappa:=\limsup_{n\to\infty}\frac{\log|I_n|}n
=\Re f_0(\mu)=\log\frac{2^{2(a+b-bc)}|\mu+c|^{bc}|\mu-c|^{bc}}
{|\mu+1|^{a+b}|\mu-1|^{a+b}},
\tag21
$$
gde $\mu$~-- vewestvenny\char"1A{} koren\char"7E~$\mu_1$
mnogoqlena~\thetag{18} iz intervala $(c,+\infty)$
v sluqae $b=1$ i koren\char"7E{} e1togo mnogoqlena
v oblasti~$\Im\tau>0$ s maksimal\char"7Eno
vozmozhno\char"1A{} qast\char"7Eyu $\Re\mu$
v sluqae $b>1$. Dlya $b=1$ verhni\char"1A{}
predel v~\thetag{21} mozhno zamenit\char"7E{}
na obyqny\char"1A.
\endproclaim

\demo{\cyrc\indent Dokazatel\char"7Estvo}\cyrm
Vse rexeniya uravneniya~\thetag{17} v oblasti $\Re\tau>0$
dlya neqetnyh $\lambda=k$ odnovremenno yavlyayut\hbox{s}ya kornyami
mnogoqlena~\thetag{18}. Rutinnaya proverka pokazyvaet, qto
uslovie~\thetag{19} obespeqivaet vozrastanie funkcii
$\Re f_0(\tau)$ kak funkcii ot~$\Re\tau$ (ili, qto to zhe samoe,
ot~$\lambda$) na krivo\char"1A~\thetag{20} v oblasti $\Re\tau>0$,
$\Im\tau\ge0$; poe1tomu na asimptotiku integrala~\thetag{15}
vliyayut tol\char"7Eko $J_{n,\pm1}$ v sluqae $b=1$ i $J_{n,\pm(b-2)}$
v sluqae $b>1$. Primenenie lemm~5 i~3 privodit k trebuemomu
sootnoxeniyu~\thetag{21}.
\enddemo

\subhead\cyrb\indent
5. Ocenki koe1fficientov line\char"1Anyh form
\endsubhead\cyrm
Dlya veliqin~$B_{k,j}$, $k=0,\pm1,\allowmathbreak\dots,\pm n$,
$j=0,1,\dots,a-1$, opredelennyh v~\thetag{9},
spravedlivy ocenki
$$
|B_{k,j}|
\le\bigl(2(a+bc-b)n\bigr)^j\cdot\max_{k=0,\pm1,\dots,\pm n}|B_{k,0}|
=\bigl(2(a+bc-b)n\bigr)^j\cdot\frac{(cn)!^{2b}(2n)!^{a+b-bc}}{n!^{2(a+b)}}.
$$
Pol\char"7Ezuyas\char"7E{} sootnoxeniyami~\thetag{10} i formulo\char"1A{} Stirlinga,
poluqaem sleduyuwee ut\-verzhdenie.

\proclaim{\cyrc\indent Lemma 7}\cyri
Dlya koe1fficientov $A_s=A_{s,n}$ line\char"1Anyh
form~\thetag{5} spravedliva ocenka
$$
\gathered
\varlimsup_{n\to\infty}\frac{\log|A_{s,n}|}n
\le2bc\log c+2(a+b-bc)\log2,
\\
\text{\cyri $s=0$ ili $s=b+1,\dots,a+b-1$ neqetno}.
\endgathered
$$
\endproclaim\cyrm

Kak neslozhno pokazat\char"7E, integraly
$$
\frac1{\pi i}\int_{iM-\infty}^{iM+\infty}
\biggl(\frac{\sin\pi t}\pi\biggr)^kR(t)\,\d t,
\qquad k=2,4,6,\dots,a-1,
\tag22
$$
vdol\char"7E{} gorizontal\char"7Eno\char"1A{} pryamo\char"1A{}
$\Im t=M$, gde $M>0$ proizvol\char"7Eno,
yavlyayut\hbox{s}ya line\char"1Anymi kombinaciyami
koe1fficientov $A_{b+2},\dots,A_{a+b-1}$
form~\thetag{5}. Poe1to\-mu primenenie asimptotiki gamma-funkcii
v oblasti $\Im t\ge M_0>0$ (sm.~\cite{B, \S\,6.5})
i metoda perevala k integralam~\thetag{22}
utoqnyaet (neznaqitel\char"7Eno) ocenku lemmy~7.

\proclaim{\cyrc\indent Lemma 8}\cyri
Pust\char"7E{} vewestvenny\char"1A{} koren\char"7E{}
$\mu_1\in(c,+\infty)$ mnogoqlena~\thetag{18}
udov\-letvoryaet usloviyu~\thetag{19} i $\eta\in(0,+i\infty)$~--
minimal\char"7Eny\char"1A{} po absolyutno\char"1A{}
veliqine mnimy\char"1A{} koren\char"7E{} e1togo mnogoqlena.
Togda dlya koe1fficientov $A_s=A_{s,n}$ line\char"1A\-nyh
form~\thetag{5} spravedliva ocenka
$$
\gathered
\varlimsup_{n\to\infty}\frac{\log|A_{s,n}|}n
\le\Re f_0(\eta)
=\log\frac{2^{2(a+b-bc)}|\eta+c|^{bc}|\eta-c|^{bc}}
{|\eta+1|^{a+b}|\eta-1|^{a+b}},
\\
\text{\cyri $s=0$ ili $s=b+1,\dots,a+b-1$ neqetno},
\endgathered
$$
priqem v sluqae $s=a+b-1$ verhni\char"1A{} predel
mozhno zamenit\char"7E{} na obyqny\char"1A{}
i neravenstvo prevrawaet\hbox{s}ya v ravenstvo.
\endproclaim

\subhead\cyrb\indent
6. Dokazatel\char"7Estvo osnovnyh rezul\char"7Etatov
\endsubhead\cyrm
Soglasno lemmam~1,~6 v slu\-qae $-b\varpi_c+2(a+b-1)+\varkappa<0$
sredi qisel $\zeta(s)$, gde $s$~neqetno i $b<s<a+b$,
imeet\hbox{s}ya po kra\char"1Ane\char"1A{} mere odno irracional\char"7Enoe. Vybiraya
$a=19$, $b=3$, $c=3$;
$a=33$, $b=5$, $c=3$ i $a=47$, $b=7$, $c=3$
sootvet\hbox{s}tvenno dlya kazhdogo iz naborov v~\thetag{1},
poluqaem teoremu~1.
V teoreme~2 dlya kazhdogo neqetnogo~$b\ge1$ polagaem $a=7b$, $c=3$;
pri e1tom koren\char"7E{} mnogoqlena~\thetag{18}
na intervale $(3,+\infty)$ sovpadaet s kornem
$\mu_1\approx3.02472$ mnogoqlena
$(\tau+3)(\tau-1)^8-(\tau-3)(\tau+1)^8$
i
$\kappa+2(a+b-1)-b\varpi_c
<\Re f_0(\mu_1)+16b-b\varpi_3
<-0.047\cdot b<0$.

V sluqae $b=1$ kriteri\char"1A{} line\char"1Ano\char"1A{} nezavisimosti iz~\cite{N1}
tak zhe, kak i v~\cite{R1}, pozvolyaet ocenit\char"7E{}
veliqinu~$\delta(a)$ snizu:
$$
\delta(a)\ge1-\frac{\kappa(a,c)+2a-\varpi_c}
{2c\log c+2(a-c+1)\log2+2a-\varpi_c},
\tag23
$$
gde veliqina $\kappa=\kappa(a,c)$ zadaet\hbox{s}ya sootnoxeniem~\thetag{21}.
Polagaya $a=145$, $c=21$ i $a=1971$, $c=131$, soglasno~\thetag{23}
poluqaem ocenki $\delta(145)\ge3$, $\delta(1971)\ge4$; krome togo,
$\delta(3)=2$ vvidu~\cite{A}. E1to dokazyvaet kak teoremu~3,
tak i teoremu~4 dlya $a<24999$. Dlya neqetnyh $a\ge20737=12^4+1$
my dokazyvaem bolee sil\char"7Enuyu, qem~\thetag{2}, ocenku
$\delta(a)>\log_{12}a$ ili, qto to zhe samoe,
$$
\delta(12^m+1)>m, \qquad m=4,5,6,\dots,
\tag24
$$
vybiraya $c=2\cdot\[a/(3m^2)\]+1$ dlya kazhdogo $a=12^m+1$
v~\thetag{23}. Pri $m=4,5,6,7$ ocenka~\thetag{24}
proveryaet\hbox{s}ya neposredstvenno; pri $m\ge8$ my pol\char"7Ezuemsya
trivi\-al\char"7Eno\char"1A{} ocenko\char"1A{} pravo\char"1A{} qasti v~\thetag{23}.
E1to zaverxaet dokazatel\char"7Estvo teoremy~4.


\bgroup
\nobreak\vskip12pt minus6pt\indent\eightpoint\smc
\hbox to100mm{\vbox{\hsize=100mm%
\leftline{\cyrms Moskovski\char"1A{} gosudarstvenny\char"1A{}
universitet im.~M.\,V.~Lomonosova}
\leftline{{\it E-mail\/}: {\tt wadim\@ips.ras.ru}}
}}\par
\egroup
\end